\documentclass{amsart}

\usepackage{amsfonts}
\usepackage{amsmath}
\usepackage{amsthm}
\usepackage{amssymb}
\usepackage[english]{babel}
\usepackage{graphics}
\usepackage{latexsym}
\usepackage{longtable} \usepackage{mathrsfs}
\usepackage{graphicx}
\usepackage{wrapfig} 
\usepackage[pdftex,colorlinks=true]{hyperref}

\newtheorem{theorem}{Theorem}
\newtheorem{corollary}[theorem]{Corollary}
\newtheorem{lemma}[theorem]{Lemma}

\newtheorem{claim}[theorem]{Claim}
\newtheorem{example}[theorem]{Example}

\theoremstyle{definition}
\newtheorem{definition}[theorem]{Definition}
\newtheorem{remark}[theorem]{Remark}

\newcommand{\mL}{\mathcal{L}}
\newcommand{\mH}{\mathcal{H}}

\newcommand{\A}{\textbf{A}}
\newcommand{\R}{\mathbb{R}}

\newcommand{\X}{\textbf{X}}
\newcommand{\Y}{\textbf{Y}}
\newcommand{\Z}{\textbf{Z}}
\renewcommand{\H}{\textbf{H}}

\newcommand{\noi}{\noindent}
\newcommand{\ms}{\medskip}

\newcommand{\al}{\alpha}
\newcommand{\be}{\beta}
\newcommand{\ga}{\gamma}
\newcommand{\de}{\delta}
\newcommand{\De}{\Delta}

\newcommand{\si}{\sigma}
\newcommand{\la}{\lambda}
\newcommand{\La}{\Lambda}
\newcommand{\ka}{\kappa}
\newcommand{\Om}{\Omega}

\newcommand{\larrow}{\longrightarrow}
\newcommand{\ot}{\otimes}

\newcommand{\p}{\partial}
\newcommand{\sub}{\subseteq}

\newcommand{\by}{\times}

\newcommand{\sgn}{\textrm{sgn}}
\newcommand{\ess}{\textrm{ess}}

\newcommand{\bt}{\begin{theorem}}\newcommand{\et}{\end{theorem}}
\newcommand{\bd}{\begin{definition}}\newcommand{\ed}{\end{definition}}
\newcommand{\bl}{\begin{lemma}}\newcommand{\el}{\end{lemma}}
\newcommand{\beq}{\begin{equation}}\newcommand{\eeq}{\end{equation}}
\newcommand{\bc}{\begin{claim}}\newcommand{\ec}{\end{claim}}
\newcommand{\bex}{\begin{example}}\newcommand{\eex}{\end{example}}
\newcommand{\bcor}{\begin{corollary}}\newcommand{\ecor}{\end{corollary}}
\newcommand{\bp}{\begin{proof}}\newcommand{\ep}{\end{proof}}

\newcommand{\BPL}{\medskip \noindent \textbf{Proof of Lemma} }

\newcommand{\BPT}{\medskip \noindent \textbf{Proof of Theorem} }

\numberwithin{equation}{section}

\begin{document}

\title[Dirichlet Problem for Fully Nonlinear Systems]{On the Dirichlet Problem for Fully Nonlinear Elliptic Hessian Systems}

\author{Nikos Katzourakis}
\address{Department of Mathematics and Statistics, University of Reading, Whiteknights, PO Box 220, Reading RG6 6AX, Berkshire, UK}
\email{n.katzourakis@reading.ac.uk}


\date{}



\begin{abstract} We consider the problem of existence and uniqueness of strong solutions $u: \Omega \subset \mathbb{R}^n \longrightarrow \mathbb{R}^N$ in $(H^{2}\cap H^{1}_0)(\Omega)^N$ to the problem
\[\label{1} \tag{1}
\left\{
\begin{array}{l}
F(\cdot,D^2u ) \,=\,  f, \ \ \text{ in }\Omega,\\
\hspace{31pt} u\,=\, 0, \ \  \text{  on }\partial \Omega,
\end{array}
\right.
\]
when $ f\in L^2(\Omega)^N$, $F$ is a Carath\'eodory map and $\Omega$ is convex. \eqref{1} has been considered by several authors, firstly by Campanato and under Campanato's ellipticity condition. By employing a new weaker notion of ellipticity introduced in recent work of the author \cite{K2} for the respective global problem on $\R^n$, we prove well-posedness of \eqref{1}. Our result extends existing ones under hypotheses weaker than those known previously. An essential part of our analysis in an extension of the classical Miranda-Talenti inequality to the vector case of 2nd order linear hessian systems with rank-one convex coefficients.
\end{abstract}

\maketitle

\section{Introduction} \label{section1}

Let $\Om \sub \R^n$ be an open bounded $C^2$ convex set, $n,N\geq 2$. Let also
\[
F\ :\ \Om \by \R^{Nn^2}_s  \larrow  \R^N
\]
be a Carath\'eodory map, namely $x\mapsto F(x,\X) $  is measurable, for every $\X \in \R^{Nn^2}_s$ and $\X\mapsto F(x,\X)$ is continuous, for almost every $x\in \Om \sub \R^{n}$. 

In this paper we consider the problem of existence and uniqueness of strong a.e.\ solutions $u : \Om \sub \R^n \larrow \R^N$ in $(H^{2}\cap H^{1}_0)(\Omega)^N$ to the following Dirichlet problem:
\beq  \label{1.1}
\left\{
\begin{array}{l}
F(\cdot,D^2u ) \,=\,  f, \ \ \text{ in }\Om,\\
\hspace{36pt} u\,=\, 0, \ \   \text{  on }\p \Om,
\end{array}
\right.
\eeq
when $ f\in L^2(\Om)^N$. In the above, $D^2u(x) \in \R^{Nn^2}_s$ is the hessian tensor of $u$ at $x$ and $Du(x) \in \R^{Nn}$ is the gradient matrix. In the sequel we will employ the summation convention in repeated indices when $i,j,k,...$ run in $\{1,...,n\}$ and $\al,\be,\ga,...$ run in $\{1,...,N\}$, while $\R^{Nn^2}_s$ is the vector space $\{\X \in \R^{Nn^2} : \X_{\al ij}=\X_{\al ji}\}$ into which the hessians of our maps are valued. The standard bases of $\R^n$, $\R^N$, $\R^{Nn}$ and $\R^{Nn^2}_s$ will be denoted by $\{e^i\}$, $\{e^\al\}$, $\{e^{\al} \ot e^i \}$ and $\{e^{\al} \ot e^i \ot e^j\}$ respectively, ``$\ot$" denotes the tensor product and we will write
\[
\text{$x=x_ie^i$, \ $u=u_\al e^\al$, \ $Du=(D_iu_\al) \, e^{\al} \ot e^i$, \ $D^2u=(D^2_{ij}u_\al)\, e^{\al}\ot e^i \ot e^j$.}
\]
Moreover, all the norms ``$|\cdot|$" appearing will always be the euclidean, e.g.\ on $\R^{Nn^2}_s$ we use $|\X|^2=\X :\X$ etc.

The problem \eqref{1.1} has been considered before by several authors and with different degrees of generality. The first one to address it was Campanato \cite{C1}-\cite{C4} under a strong ellipticity condition which we recall later. Subsequent contributions to this problem and problems relevant to Campanato's work can be found in Tarsia \cite{Ta1}-\cite{Ta5}, Fattorusso-Tarsia \cite{FT1}-\cite{FT4}, Buica-Domokos \cite{BD}, Domokos \cite{D}, Palagachev \cite{Pa1,Pa2}, Palagachev-Recke-Softova \cite{PRS}, Softova \cite{S} and Leonardi \cite{Le}. However, all vectorial contributions, even the most recent ones \cite{FT1,FT2} (wherein they consider systems of the form $F(\cdot,u,Du,D^2u)=f$) are based on Campanato's original restrictive ellipticity notion, or a small extension of it due to Tarsia \cite{Ta5}.

The main consequence of Campanato's ellipticity is that the nonlinear operator $F[u]:=F(\cdot,D^2u)$ is ``near" the Laplacian $\De u$. Nearness is a functional analytic notion also introduced by Campanato in order to solve the problem, which roughly says that operators near those with ``good properties" like bijectivity inherit these properties. In the case at hand, nearness implies unique solvability of \eqref{1.1} in $(H^{2}\cap H^{1}_0)(\Om)^N$, by the unique solvability of the Poisson equation $\De u =f$ in $(H^{2}\cap H^{1}_0)(\Om)^N$ and a fixed point argument. Campanato's ellipticity relates to the Cordes condition (see Cordes \cite{Co1,Co2} and also Talenti \cite{T} and Landis \cite{L}). 

Although Campanato's condition is stringent, it should be emphasised that in general it is not possible to obtain solvability in the class of strong solutions with the mere assumption of uniform ellipticity. Well-known counterexamples which are valid even in the linear scalar case of the second order elliptic equation 
\[
A_{ij}(x)D^2_{ij}u(x)\, =\, f(x) 
\]
with $A_{ij} \in L^\infty(\Om)$ imply that the standard uniform ellipticity $ A\geq \nu I$ does not suffice to guarantee well posedness of the Dirichlet problem when $n>2$ and more restrictive conditions are required (see e.g.\ Ladyzhenskaya-Uraltseva \cite{LU}).

In this work we prove well posedness of \eqref{1.1} in the space $(H^{2}\cap H^{1}_0)(\Om)^N$ for any $f\in L^2(\Om)^N$ under a \emph{new ellipticity condition on $F$} which is strictly weaker than the Campanato-Tarsia notion. This new notion has been introduced in the very recent paper of the author \cite{K2} in order to study the case of the global problem on the whole space $\R^n$ for the same fully nonlinear hessian system:
\[
F(\cdot, D^2u)\, =\, f, \ \ \ u\, :\, \R^n \larrow \R^N.
\]
The relevant first order global problem $F(\cdot, Du)\, =\, f$ has also been studied in \cite{K1}, which is a non-trivial generalisation of the Cauchy-Riemann equations. The idea of our weaker notion is to require $F$ to be ``near" a general 2nd order elliptic system with constant coefficients which satisfies the Legendre-Hadamard condition, instead of being ``near" the Laplacian.

More precisely, our starting point for the system $F(\cdot,D^2u)=f$ is based on the analysis of the simpler case of $F$ linear in $\X$ and independent of $x$, that is when
\beq \label{1.2}
F_\al(x,\X)\, =\, \A_{\al i \be j} \X_{\be i j}.
\eeq
Here $\A$ is a linear symmetric operator $\A : \R^{N n}\larrow \R^{Nn}$:
\[
\A \in \R^{Nn \by Nn}_s, \ \  \text{ i.e. }\ \A_{\al i \be j}\, =\, \A_{\be j \al i}.
\]
For $F$ as in \eqref{1.2}, the system $F(\cdot,D^2u)=f$ becomes
\[
\A_{\al i \be j} D^2_{ij}u_\be \,=\, f_\al.
\]
By introducing the contraction operation $\A :\Z := (\A_{\al i \be j}\Z_{\al ij})e^\al$ (which extends the trace inner product $\Z:\Z=\Z_{\al i j}\Z_{\al i j}$ of $\R^{Nn^2}_s$), we will write it compactly as
\beq \label{1.3}
\A: D^2u \,=\, f.
\eeq
The appropriate notion of ellipticity in this case is that the quadratic form arising from the operator $\A$
\beq \label{1.4}
\begin{array}{c}
\A \ :\ \ \R^{Nn} \by \R^{Nn}\ \larrow \R,\ms\\ 
\A: P\ot Q \, := \, \A_{\al i \be j}P_{\al i}Q_{\be j},
\end{array}
\eeq
is (strictly) rank-one convex on $\R^{Nn}$, that is
\beq \label{1.5}
\A : \eta \ot a \ot \eta \ot a \, \geq \, \nu |\eta|^2|a|^2,
\eeq
for some $\nu >0$ and all $\eta \in \R^N, \ a\in \R^n$. For brevity, we will say \emph{``$\A$ is rank-one positive"} as a shorthand of the statement \emph{``the symmetric quadratic form defined by $\A$ on $\R^{Nn}$ is rank-one convex"}. Our ellipticity assumption for general $F$ is given in the following:

\begin{definition}[K-Condition] \label{def2} \textit{ Let $\Om \sub \R^n$ be open and $F :  \Om \by \R^{Nn^2}_s \larrow  \R^N$ a Carath\'eodory map. We say  that $F$ is elliptic (or that the PDE system $F(\cdot,D^2u)=f$ is elliptic) when there exist
\[
\begin{array}{c}
\A \in \R^{Nn \by  Nn}_s\ \text{ rank-one positive},\smallskip\\
\al \in L^\infty(\Om), \ \al>0 \text{ a.e.\ on }\Om, \ 1/\al \in L^\infty(\Om),\smallskip\\
\be,\ga>0\ \text{ with }\be+\ga<1,
\end{array} 
\]
such that
\beq \label{2.3}
\left| \A:\Z\, -\, \al(x)\Big(F(x,\X+\Z) -F(x,\X)\Big) \right| \ \leq\,  \be\, \nu(\A)|\Z|\ + \ \ga\, |\A:\Z|,
\eeq
for all $\X, \Z \in \R^{Nn^2}_s$ and a.e.\ $x\in \Om$.}
\end{definition}

In the above definition $\nu(\A)$ is the ellipticity constant of $\A$:
\beq \label{1.7}
\nu(\A)\, :=\, \min_{|\eta|=|a|=1} \big\{\A : \eta \ot a \ot \eta \ot a \big\}. 
\eeq
By taking 
as $\A$ the monotone tensor
\[
\A_{\al i \be j}\ =\ \de_{\al \be}\de_{ij}, 
\]
we reduce to a condition equivalent to Tarsia's notion, and by further taking $\al(x)$ constant we reduce to Campanato's notion:
\beq \label{1.8}
\left| \Z:I\, -\, \al \Big(F(x,\X+\Z) -F(x,\X)\Big) \right| \ \leq\,  \be |\Z|\ + \ \ga|\Z:I|,
\eeq
In \eqref{1.8} we have used the obvious contraction operation $\X:X:=(\X_{\al i j}X_{ij})e^\al$. Our new ellipticity notion \eqref{2.3} relaxes \eqref{1.8} substantially: a large class of nonlinear operators which are elliptic are of the form
\[
F(x,\X)\, :=\, g^2(x)\A:\X \ + \ G(x,\X)
\]
where $\A$ is rank-one positive, $g,1/g \in L^\infty(\Om)$ and $G$ is \emph{any nonlinear} map, measurable with respect to the first argument and Lipschitz with respect to the second argument, with Lipschitz constant of $G(x,\cdot)/g^2(x)$ smaller than $\nu(\A)$ (see Example 5 in \cite{K2}). In particular, any $F \in C^{1}\big(\R^{Nn^2}_s\big)^N$ such that $F'(0)$ is rank-one positive and the Lipschitz constant of $\X\mapsto F(\X)-F'(0):\X$ is smaller than $\nu(F'(0))$, is elliptic in the sense of Definition \ref{def2}. On the other hand, even if $F$ is linear, $F(\X)=\A:\X$ and in addition $\A$ defines a \emph{strictly convex} quadratic form on $\R^{N\by n}$, that is when
\[
\A : Q \ot Q \, \geq c^2 |Q|^2,\ \ \ Q\in \R^{N\by n},
\]
then $F$ may \emph{not be elliptic} in the Campanato-Tarsia sense (see Example 6 in \cite{K2}).

The program we deploy herein is the following: we first solve the Dirchlet problem \eqref{1.1} in the linear case with constant coefficients for $F(\X)=\A:\X$. This is a simple application of classical variational and regularity results and is recalled in Section \ref{section2}. Next, in Section \ref{section3} we establish a crucial ingredient of our analysis: a sharp estimate in $(H^2\cap H^1_0)(\Om)^N$ for linear hessian operators with rank-one positive constant coefficients which is an extension of the classical Miranda-Talenti inequality (see \cite{M}, \cite{T} and also \cite{DP}). Namely, in Lemma \ref{pr2} we show that if $\A \in \R^{Nn \by Nn}_s$ is rank-one positive and $\Om$ is a bounded convex $C^2$ domain, then for all maps $u :\Om \sub \R^n \larrow \R^N$ in $(H^2\cap H^1_0)(\Om)^N$, we have the estimate
\beq \label{1.12}
\big\|D^2u \big\|_{L^2(\Om)} \, \leq \, \frac{1}{\nu(\A)}\big\|\A:D^2u \big\|_{L^2(\Om)} .
\eeq
The inequality \eqref{1.12} is a vectorial non-monotone extension of the Miranda-Talenti inequality beyond the scalar case of the Laplacian of the classical result and appears to be anew result even in the scalar case. By choosing $N=1$ and as $\A$ the identity of $\R^n$
\[
\A_{\al i\be j}\, =\, \de_{ij}
\]
we reduce to the classical scalar case with $\A:D^2 u=\De u$ and $\nu(\A)=1$. However, we point out a weakness of our result: we were able to prove \eqref{1.12} only under an extra technical assumption on the minors of the 4th order tensor $\A$, whose necessity unfortunately we have not been able to verify. This extra assumption trivialises in the \emph{scalar} case; indeed, \eqref{1.12} holds when $N=1$ for \emph{any} positive matrix $\A\in \R^{n^2}_s$ without further restrictions. Notwithstanding, even under the extra condition, \eqref{1.12} is still a genuine extension to a new realm. The proof builds on the Miranda-Talenti identity
\beq \label{1.13}
\int_\Om \Big\{\big| D^2v \big|^2\,-\ (\De v)^2 \Big\}d\mL^n \, =\, (n-1)\int_{\p \Om} |Dv|^2\,\H \cdot N\, d\mH^{n-1}
\eeq
valid for scalar functions $v\in (H^2\cap H^1_0)(\Om)$, where $N$ is the outwards pointing unit vector field of $\p\Om$ and $\H$ is the mean curvature vector. The proof of \eqref{1.13} is recalled in the Appendix, for the sake of completeness, for the convenience of the reader and also because with the exception of the book \cite{MPS}, this important result is not easy to find in English-written literature.

Next, in Section \ref{section4} we consider the general case of fully nonlinear $F$ satisfying Definition \ref{def1}  (Theorem \ref{th2}). The idea is to use our ellipticity notion which serves as a ``perturbation device" and employ Campanato's theorem of bijectivity of near operators in order to connect the nonlinear to the linear problem.  Campanato's result is taken from \cite{C5} and the short proof is recalled in the Appendix for the convenience of the reader (Theorem \ref{th3}). Our analysis follows very similar lines to the respective proof of \cite{K2} for the global problem on $\R^n$, but we chose to give all the details here too. A byproduct of our method is a strong uniqueness estimate in the form of a \emph{comparison  principle} for the distance of any solutions in terms of the distance of the right hand sides of the equations. Moreover, in Section \ref{section5} we discuss a result of stability type for the Dirichlet problem over bounded domains, along the lines of respective result of \cite{K2} for global solutions.

We note that Campanato's notion of nearness has been relaxed by Buica-Domokos in \cite{BD} to a ``weak nearness", which still retains most of the features of (strong) nearness. In the same paper, the authors also use an idea similar to ours, namely a fully nonlinear operator being ``near" a general linear operator, but they implement this idea only in the scalar case. 

We conclude this introduction by noting that \eqref{1.1} has been studied also when $F$ is \emph{coercive instead of elliptic}. By using the analytic Baire category method of Dacorogna-Marcellini \cite{DM}, one can prove that, under certain structural and compatibility assumptions, the Dirichlet problem has \emph{infinitely many} strong a.e.\ solutions in $W^{2,\infty}(\Om)^N$. This method is the ``geometric counterpart" of Gromov's Convex Integration. However, ellipticity and coercivity of $F$ are, roughly speaking, mutually exclusive and this method does not in general give uniqueness. On the other hand, the bibliography on the scalar theory of elliptic equations is vast, for both classical/strong a.e.\ solutions, (see e.g.\ Gilbarg-Trudinger \cite{GT}) as well as for viscosity solutions of degenerate equations  (for an elementary intorduction see e.g.\ \cite{K}). However, except for the broad theory for divergence structure systems (see e.g.\ \cite{GM}), for fully nonlinear systems the existing theory is very limited and this applies even to linear non-variational systems.

\section{Preliminaries and well posedness of the linear problem} \label{section2}

We begin by recalling that in \cite{K2} we used as primary ellipticity notion that of Definition \ref{def1} below. The latter is equivalent to the $K$-condition of Definition \ref{def2}, under the assumption of global Lipschitz continuity of $F$ with respect to the second argument.

\begin{definition}[Ellipticity, cf.\ \cite{K2}] \label{def1} Let $F : \Om \by \R^{Nn^2}_s \larrow  \R^N$ be a Carath\'eodory map. We call $F$ (or the system $F(\cdot,D^2u)=f$) elliptic when there exist 
\[
\begin{array}{c}
\A\,  \in \, \R^{Nn \by Nn}_s, \text{ rank-one positive},\\
\la>\ka>0,\\
\text{$\al \in L^\infty(\Om)$, $\al>0$ a.e.\ on $\Om$ and $1/\al \in L^\infty(\Om)$,}
\end{array}
\]
such that
\beq \label{1.6}
(\A :\Z)^\top\Big[ F(x,\X+\Z) - F(x,\X)\Big] \ \geq\ \frac{\la}{\al(x)}|\A:\Z|^2\, -\ \frac{\ka}{\al(x)}\nu(\A)^2|\Z|^2,
\eeq
for all $\X,\Z \in \R^{Nn^2}_s$ and a.e.\ $x\in \Om \sub \R^n$. 
\end{definition}

\noi We recall that $\nu(\A)$ is the ellipticity constant of $\A$ and is given by \eqref{1.7}. The inequality \eqref{1.6} can be seen as a pseudo-monotonicity condition (see also \cite{Ta3, Ta4}).

\ms

We now consider the question of unique solvability of the Dirichlet problem \eqref{1.1} in the case of linear systems
\[
\left\{
\begin{array}{l}
\A:D^2u \,=\,  f, \ \ \text{ in }\Om,\\
\hspace{30pt} u\,=\, 0, \ \   \text{  on }\p \Om,
\end{array}
\right.
\]
for any $ f\in L^2(\Om)^N$, when $\A\in \R^{Nn \by Nn}_s$ is strictly rank-one positive. This is a standard application of the direct method of Calculus of Variations (see e.g.\  Dacorogna \cite{D}) in order to get existence of a weak solution of the Euler-Lagrange equation
\[
-D_i \Big( \A_{\al i \be j}D_ju_\be \Big)\, +\, f_a\, =\ 0,
\]
by minimising the functional
\[
E(u,\Om)\, =\, \int_\Om \Big( \A:Du(x) \ot Du(x) \,+\, f(x)\cdot u(x)\Big)\, dx
\]
in $H^1_0(\Om)^N$ and then apply regularity theory.

\begin{lemma}[Well posedness of the linear problem] \label{le3} Let $n,N\geq 2$ and $\Om \sub \R^n$ a bounded domain with $C^2$ boundary. Then, for any $\A\in \R^{Nn \by Nn}_s$ strictly rank-one positive and any  $ f\in L^2(\Om)^N$, the problem
\[
\left\{
\begin{array}{l}
\A:D^2u \,=\,  f, \ \ \text{ in }\Om,\\
\hspace{30pt} u\,=\, 0, \ \   \text{  on }\p \Om,
\end{array}
\right.
\]
has a unique strong solution in the space $(H^2\cap H^1_0)(\Om)^N$, which solves the PDE system a.e.\ on $\Om$. Moreover, the solution $u$ satisfies the estimate
\[
\|u\|_{L^2(\Om)}\, +\,\|Du\|_{L^2(\Om)}\, +\,\|D^2u\|_{L^2(\Om)}\ \leq \ C\|f\|_{L^2(\Om)},
\]
with $C>0$ depending only on $\Om'$, $\Om$ and $\A$.
\end{lemma}

\BPL \ref{le3}. The proof can be found e.g.\ in Giaquinta-Martinazzi \cite{GM} spread over the pages 55-72.      \qed

\ms

\begin{remark}[Equivalent norms on $(H^2\cap H^1_0)(\Om)^N$] \label{rem} We record the standard fact that Poincar\'e inequality in $H^1_0$ and interpolation inequalities in $L^2$ (see e.g.\ Gilbarg-Trudinger \cite{GT}) imply that two equivalent norms on $(H^2\cap H^1_0)(\Om)^N$ are
\[
\|D^2u\|_{L^2(\Om)}\ \approx\   \|u\|_{H^2(\Om)}\, :=\, \|u\|_{L^2(\Om)}\, +\, \|Du\|_{L^2(\Om)}\, +\, \|D^2u\|_{L^2(\Om)}. 
\]
\end{remark} 

\ms

\section{The generalised Miranda-Talenti inequality for elliptic systems with constant coeffients}
\label{section3}

In this section we establish the estimate \eqref{1.12} in Lemma \ref{pr2} below. This is an extension of the Miranda-Talenti inequality from the case of the Laplacian to the case of general $\A$. We note that this result, even in the classical case of the Laplacian, is \emph{non-trivial}. \emph{The fact that we do not restrict the gradient to vanish on the boundary is an essential difficulty}. The inequality \eqref{1.12} in the smaller space 
\[
H^2_0(\Om)^N \, =\, \overline{\, C^\infty_c(\Om)^N\, }^{\|\cdot\|_{H^2}}
\]
(instead of $(H^2\cap H^1_0)(\Om)^N$) does not require boundary regularity and holds for rank-one positive $\A$ (that is $\nu(\A)>0$) without extra conditions. The proof in $H^2_0(\Om)^N$ follows by applying the Fourier transform, Plancherel's theorem, the properties of rank-one convexity and an approximation argument. This is done in \cite{K2}, althought the result is stated directly for the whole space $\R^n$. On the other hand, standard global $L^2$ regularity theory for linear systems (see \cite{GM}) says that \eqref{1.12} is always true in $(H^2\cap H^1_0)(\Om)^N$ for any $C^2$ domain $\Om$ without curvature restrictions, \emph{but with a perhaps larger universal constant, instead of the sharp value $1/\nu(\A)$}.

As we have already pointed out in the introduction, in the vectorial case we need an extra technical condition on $\A$ except for rank-one convexity. This restriction is void in the scalar case $N=1$ and the scalar version of \eqref{1.12} is true for a general positive $\A \in \R^{n \by n}_s$. The assumption we need is the next

\ms

\noi \textbf{Structural Hypothesis (SH).} \textit{Let $n,N\geq 1$. Consider a tensor $\A \in \R^{Nn \by Nn}_s$, which we view as a linear map
\[
\A\ :\ \R^{n^2} \larrow \R^{N^2},\ \ \ \ \A :X \, =\, \big(\A_{\al i \be j}X_{ij}\big)e^{\al} \ot e^\be. 
\]
We assume there exist matrices $B^1,...,B^N$ in $\R^{N^2}_s$ and $A^1,...,A^N$ in $\R^{n^2}_s$ such that $\A$ can be written as
\[
\A \ =\ B^1\ot A^1\ + \ \cdots\  +\  B^N\ot A^N
\]
and $B^\ga, A^\ga$ satisfy
\[
\left\{
\begin{array}{c}
R(B^\ga) \,\bot\,  R(B^\de), \ \ \text{ for }\ga \neq \de, \ms\\
B^1\, + \ \cdots\  +\, B^N >\, 0, \ms\\ 
A^1 >\,0, \ \ldots\ , A^N >\,0, \ms \\
\dim \left( \bigcap_{\ga=1}^N N\big(A^\ga - \la_1(A^\ga)I \big)\right) \, \geq\, 1.
\end{array}
\right.
\]
In the above, $\la_1$ denotes the smallest eigenvalue, $R$ denotes the range and $N$ denotes the nullspace.
} 

\begin{remark}[(SH) $\Longrightarrow$ rank-one convexity] Every tensor $\A$ which satisfies (SH) is necessarily rank one positive: indeed, for any $\eta \in \R^N$ and $a\in \R^n$, we have
\[
\begin{split}
\A :\eta \ot a \ot \eta \ot a\, &=\, \big( B^\ga : \eta \ot \eta \big) \big( A^\ga : a \ot a \big)\\
&\geq\, \Big[ \big(B^1+ \cdots+B^N \big) : \eta \ot \eta \Big] \min_{\ga=1,...,N}  \big\{ A^\ga : a \ot a \big\}\\
&\geq \, \Big\{\la_1\big(B^1+ \cdots+B^N\big) \min_{\ga=1,...,N}\la_1(A^\ga) \Big\} \, |\eta|^2|a|^2
\end{split}
\]
and the quantity in the last bracket is strictly positive by the positivity of the sum $B^1+\cdots+B^N$ and of $A^1$, ..., $A^N$. Let us also record the identity 
\[
\nu(A)\, =\, \la_1(A)\, =\, \min_{|a|=1} \big\{ A:a \ot a \big\},
\]
which is valid for any postive $A\in \R^{n^2}_s$ and we have just used it in the last step of the previous inequality. 

We also note that (SH) implies that each $B^\ga$ is non-negative, but in general may have non-trivial nullspace. However, the sum $B^1+\cdots+B^N$ is strictly positive and the direct orthogonal sum of their ranges spans the space $\R^N$.
\end{remark}

\begin{remark} The existence of plenty of on-trivial examples of $\A$'s satisfying (SH) is fairly obvious. The special case of the \emph{monotone} operator $\A = I \ot A$ with $A\in \R^{n^2}_s$ positive, that is when
\[
\A_{\al i \be j} \, =\, \de_{\al \be} A_{ij},
\]
automatically satisifes (SH) with
\[
B^1=I,\ \ B^2=\cdots=B^N=\textbf{0}, \ \ \ A^1=\cdots=A^N=A.
\]
If in addition $A=I$, that is if $\A_{\al i \be j} \, =\, \de_{\al \be} \de_{ij}$, then $\A$ gives rise to the vectorial Laplacian operator: $\A_{\al i \be j} D^2_{ij}u_\be = D^2_{ii} u_\al=\De u_\al$.
\end{remark}

\begin{lemma}[The generalised Miranda-Talenti inequality for linear hessian systems] \label{pr2}

Let $\A\in \R^{Nn \by Nn}_s$ be rank-one positive with ellipticity constant $\nu(\A)$ given by \eqref{1.7} and satisfying the structural hypothesis (SH). 

Let also $\Om\sub \R^n$ be open, convex and bounded. Then, we have the estimate
\beq \label{3.1}
\big\|D^2u \big\|_{L^2(\Om)} \, \leq \, \frac{1}{\nu(\A)}\big\|\A:D^2u \big\|_{L^2(\Om)} ,
\eeq
valid for all maps $u :\Om\sub \R^n \larrow \R^N$ in $(H^2\cap H^1_0)(\Om)^N$.
\end{lemma}

The proof of Lemma \ref{pr2} is based on the Miranda-Talenti identity (\cite{M,T})
\beq \label{3.2}
\int_\Om \Big\{\big| D^2v \big|^2\,-\ (\De v)^2 \Big\}d\mL^n \, =\, (n-1)\int_{\p \Om} |Dv|^2\,\H \cdot N\, d\mH^{n-1}
\eeq
valid for scalar functions $v\in (H^2\cap H^1_0)(\Om)$, where $N$ is the outwards pointing unit vector field of $\p\Om$, $\H$ is the mean curvature vector, $\mL^n$ is the Lebesgue measure and $\mH^{n-1}$ is the Hausdorff measure. 

\ms

\noi\textbf{The idea of the proof.} Roughly, (SH) allows to decouple $\A$ to a sum of product subtensors which are in a certain sense orthogonal to each other. For each decomposed subtensor, we can use appropriate transformations to reduce the matrices comprising it to the product of the identity matrices on $\R^N$ and $\R^n$ respecetively. Then, we can apply the classical result \eqref{3.2} to each component of the tranformed product subtensors. By reassembling all the components back together and inverting the tranformations, we get \eqref{3.1}. The idea is simple, but the proof has some technicalities.

\ms

\begin{remark}[On the convexity assumption for $\p \Om$] It is well know that if $\Om$ is convex, then the mean curvature vector points towards the interior. When $n\geq 3$, this is strictly weaker than convexity and even non-simply connected domains may satisfy it. For example, the torus $\mathbb{T}^2\sub \R^3$ can satisfy $\H \cdot N\leq 0$ if the ratio of the radii is chosen appropriately. \emph{However}, in the general case of $\A$ we are dealing with in this work, \emph{we can not in general relax the convexity requirement} for $\Om$. This will be obvious from the proof and counterexamples are easy to demostrate, but we refrain from this task.
\end{remark}

\BPL \ref{pr2}. \textbf{Step 1}. We first prove \eqref{3.1} in the scalar case of a positive $\A =A \in \R^{n^2}_s$. Since $A>0$, by the Spectral theorem we can find an orthogonal matirx $O\in O(N,\R)$ and a positive diagonal matrix $\La \in \R^{n^2}_s$ such that
\[
\begin{array}{c}
A\, =\, K \, K^\top,\ \ \ K\, =\, O\La, \ms \ms\\
\La\, =\, \left[
\begin{array}{ccc}
\sqrt{\la_1} & & \textbf{0} \\
 & \ddots & \\
\textbf{0} & & \sqrt{\la_n} 
\end{array}
\right], \ \ \ \la_i=\la_i(A), 
\ms \ms
\\
\si(A) \, =\, \big\{\la_1 ,....,\la_n \big\},\ \ \la_i\leq \la_{i+1}.
\end{array} \ \ \ \ \ \ \
\]
Namely, the entries of $\La$ are the square roots of the eigenvalues of $A$. We now show the following algebraic inequality
\beq \label{3.3}
\big| \La \,X\, \La \big| \, \geq \, \nu(A) \, |X|^2,
\eeq
which is true for any $X \in \R^{n^2}_s$. Indeed, since $\La_{ij}=0$ for $i\neq j$ and $\La_{ii}=\sqrt{\la_i}$, we have (we will now disengage the summation convention in order to avoid confusion with the repeated indices which do not sum)
\[
\begin{split}
\big| \La \,X\, \La \big|^2 \, &=\, \big( \La \,X\, \La \big) : \big( \La \,X\, \La \big)\\
&=\, \Bigg( \sum_{i,j}\sqrt{\la_i} \, X_{ij} \,  \sqrt{\la_j} \, e^i \ot e^j \Bigg) : \left(\sum_{kl}\sqrt{\la_k} \, X_{kl} \,  \sqrt{\la_l} \, e^k \ot e^l \right) \\
&=\,  \sum_{i,j,k,l} X_{ij} \, \sqrt{\la_i} \, \sqrt{\la_j} \, X_{kl} \, \sqrt{\la_k} \, \sqrt{\la_l}   \,\de_{ik} \de_{jl}\\
&=\,  \sum_{i,j} X_{ij} \, \sqrt{\la_i} \, \sqrt{\la_j}  X_{ij} \, \sqrt{\la_i} \, \sqrt{\la_j} \\
&=\,  \sum_{i,j} (X_{ij})^2 \,  {\la_i} \,  {\la_j}.  
\end{split}
\]
Since $\la_i \geq \la_1=\nu(A)$, we obtain
\[
\begin{split}
\big| \La \,X\, \La \big|^2 \, &\geq \, (\la_1)^2 \sum_{i,j} (X_{ij})^2  \\
&= \, \nu(A)^2 |X|^2
\end{split}
\]
and hence \eqref{3.3} has been established. Next, we fix $v\in C^2(\overline{\Om})\cap C^1_0(\Om)$ and set 
\[
\tilde{\Om}\,:=\, K^{-1}\Om,\ \ \ \ \tilde{v}\ :\ \tilde{\Om} \sub \R^n \larrow \R, \ \ \tilde{v}(x)\, :=\, v(Kx).
\]
Then, for any fixed $x\in \tilde{\Om}$, we have
\[
D^2_{ij}\tilde{v}(x)\, =\, D^2_{kl}v(Kx) \, K_{ki}\, K_{lj}\, =\, K^\top_{ik}\,D^2_{kl}v(Kx) \, K_{lj}
\]
and hence
\beq \label{3.4}
\De \tilde{v}(x)\, =\, D^2_{kl}v(Kx) \, K_{ki}\, K_{li}\, =\, A: D^2v(Kx).
\eeq
Moreover, we have
\[
\begin{split}
\big| D^2 \tilde{v}(x)\big|^2\, &=\, \Big| K^\top D^2 v(Kx) K\Big|^2\\
&=\, \Big| \La \left( O^\top D^2 v(Kx)\, O \right) \La \Big|^2\\
&\overset{\eqref{3.3}}{\!\!\! \geq} \nu(A)^2 \, \left| O^\top D^2 v(Kx)\, O \right|^2.
\end{split}
\]
Hence, we get
\[
\begin{split}
\big| D^2 \tilde{v}(x)\big|^2 \, 
&\geq\, \nu(A)^2 \, \left| O^\top D^2 v(Kx)\, O \right|^2\\
& =\, \nu(A)^2 \, O^\top_{ik}\,D^2_{kl}v(Kx) \, O_{lj} \, O^\top_{ip}\,D^2_{pq}v(Kx) \, O_{qj}\\
& =\, \nu(A)^2 \, D^2_{kl}v(Kx) \,D^2_{pq}v(Kx) \, \de_{pk} \de_{ql}
\end{split}
\]
which gives
\beq \label{3.5}
\big| D^2 \tilde{v}(x)\big|^2 \, \geq\, \nu(A)^2 \, \left| D^2 v(Kx) \right|^2.
\eeq

We now claim that since $\Om$ is a $C^2$ bounded convex set, $ \tilde{\Om}=K^{-1}\Om$ is $C^2$ bounded convex too. Indeed, since 
\[
K^{-1}\,=\,(O\La)^{-1}\,=\,\La^{-1}O^\top, 
\]
we have $\tilde{\Om} = \La^{-1}(O^\top \Om)$. Since $O$ is orthogonal, $O^\top \Om$ is isometric to $\Om$ and hence convex. Let us set
\[
C\,:=\, O^\top \Om.
\]
Then, $\tilde{\Om} =\La^{-1}C$ is also a convex set. To see this, note that we can find a convex function $F\in C^2(\R^n)$ such that $\{f<0\}=C$. For example, one such function is given by 
\[
F(x)\, :=\, \inf\Big\{t>-1\ :\ x\in \big((t+1)(C-\bar{x})\big) + \bar{x} \Big\}
\]
where $\bar{x} $ is any fixed point in $\Om$. Then, we consider the function
\[
\tilde{F}(x)\, :=\, F(\La x),\ \ \ \tilde{F}\in C^2(\R^n).
\]
Then, we have
\[
D_{ij}^2\tilde{F}(x)\, =\, D^2_{pq}F(\La x)\, \La_{pi}\, \La_{qj}
\]
and hence for any $a\in \R^n$, the convexity of $F$ implies
\[
\begin{split}
D^2\tilde{F}(x) : a\ot a \, &=\, D_{ij}^2\tilde{F}(x) \,a_i\, a_j\\
& =\, D^2_{pq}F(\La x)\, \La_{pi}\, \La_{qj}\, a_i\, a_j\\
&=\, D^2 F(\La x) : \big( \La a\big) \ot  \big( \La a\big) \\
&\geq\, 0.
\end{split}
\]
Hence, $\tilde{F}$ is convex too, which means the sublevel set $\{\tilde{F} <0 \}$ is convex. Moreover,
\[
\begin{split}
\tilde{\Om} \, =\, \La^{-1}C \, & =\, \La^{-1} \big\{x\in \R^n\ :  \ F(x)<0\big\}\\
& =\,  \big\{\La^{-1}x\in \R^n\ :  \ F(x)<0\big\}\\
&=\, \left\{y\in \R^n\ :  \ F(\La y)<0\right\}  \\
& =\, \big\{y\in \R^n\ :  \ \tilde{F}(y)<0\big\}.
\end{split}
\]
Thus, $\tilde{\Om}$ is convex and the conclusion follows.

Now, since $\tilde{v} \in H^2(\tilde{\Om})\cap H^1_0(\tilde{\Om})$ and $\tilde{\Om}$ is convex, we may applying the Miranda-Talenti identity \eqref{3.2} to $\tilde{v} $ and $\tilde{\Om}$ to obtain
\[
\int_{\tilde{\Om}} \big| D^2 \tilde{v}(x)\big|^2dx \, \leq \, \int_{\tilde{\Om}} \big| \De \tilde{v}(x)\big|^2dx .
\]
Hence, by using \eqref{3.4} and  \eqref{3.5}, we obtain
\[
\nu(A)^2 \int_{\tilde{\Om}}  \, \left| D^2 v(Kx) \right|^2 \, \leq \, \int_{\tilde{\Om}} \big| A: D^2v(Kx)\big|^2 dx.
\]
We conclude by using the change of variables $y=Kx$ which sends $\tilde{\Om}$ back to $\Om$ and a standard approximation argument in the Sobolev norm. Hence, our inequality \eqref{3.1} has been established in the scalar case.

\ms

\noi \textbf{Step 2.} We now start working towards the general vector case. Hence, let $\A \in \R^{Nn \by Nn}_s$ satisfy the structural hypothesis (SH) for some matrices $\{B^1,...,B^N\} \sub \R^{N^2}_s$ and $\{A^1,...,A^N\}\sub \R^{n^2}_s$. We begin by showing that \emph{we may further assume that all the positive matrices $A^1,...,A^N$ have the same first eigenvalue}: 
\[
\la_1(A^1)\, =\, \cdots\, =\, \la_1(A^N)\, =:\, \la_1 >0.
\]
Indeed, if this is not the case, we may find positive constants $c_1,...,c_N$ such that
\[
\begin{split}
\A\, &=\, B^1 \ot A^1\, +\, \cdots \, +\, B^N \ot A^N\\
&=\, (c_1B^1) \ot \frac{A^1}{c_1}\, +\, \cdots \, +\,  (c_NB^N) \ot \frac{A^N}{c_N}\\
 &=:\, \tilde{B}^1 \ot \tilde{A}^1\, +\, \cdots \, +\, \tilde{B}^N \ot \tilde{A}^N\\
\end{split}
\]
and the rescaled families of matrices $\{\tilde{B}^1,...,\tilde{B}^N\}$ and $\{\tilde{A}^1,...,\tilde{A}^N\} $ have the same properties as $\{B^1,...,B^N\} $ and $\{A^1,...,A^N\} $, but in addition there exists an $\bar{a} \in \R^n$ with $|\bar{a}|=1$ and $\la_1>0$ such that
\beq \label{3.6}
\min_{|a|=1} \big\{\tilde{A}^\ga : a \ot a \big\}\, =\,  \tilde{A}^\ga : \bar{a} \ot \bar{a} \, =\, \la_1,
\eeq
for all $\ga=1,...,N$. The existence of such an $\bar{a}$ for all the $\ga$'s is provided by (SH): by assumption, all nullspaces $N\big(A^\ga -\la_1(A^\ga)I\big)$ for $\ga=1,...,N$ intersect at least along a common line of $\R^n$. 

Next, we show that, under the previous simplification, the ellipticity constant $\nu(\A)$ of $\A$ defined by \eqref{1.7} is also given by
\beq \label{3.7}
\nu(\A)\, =\, \left(\min_{\ga=1,...,N} \min_{|a|=1} \big\{ A^\ga : a \ot a \big\}\right) \min_{|\eta|=1} \Big\{ \big(B^1+\cdots+B^N\big): \eta \ot \eta \Big\}.
\eeq
Indeed, we have 
\[
\begin{split}
\nu(\A)\, &=\,  \min_{|a|=|\eta|=1} \Big\{  \big( B^\ga \ot A^\ga \big) :  \eta \ot \eta \ot a \ot a \Big\}\\
&=\, \min_{|a|=|\eta|=1} \Big\{ \big( A^\ga : a \ot a \big)  \big( B^\ga: \eta \ot \eta \big)\Big\}\\
&\geq\, \left(\min_{\ga=1,...,N} \min_{|a|=1} \big\{ A^\ga : a \ot a \big\}\right) \min_{|\eta|=1} \Big\{ \big(B^1+\cdots+B^N\big): \eta \ot \eta \Big\},
\end{split}
\]
and conversely, by the previous arguments (see \eqref{3.6}),
\[
\begin{split}
\nu(\A)\, &=\, \min_{|a|=|\eta|=1} \Big\{ \big( A^\ga : a \ot a \big)  \big( B^\ga: \eta \ot \eta \big)\Big\}\\
&\leq\, \min_{|\eta|=1} \Big\{ \big( A^\ga : \bar{a} \ot \bar{a} \big)  \big( B^\ga: \eta \ot \eta \big)\Big\}\\
&=\, \min_{|\eta|=1} \Big\{ \la_ 1  \Big[ \big(B^1+\cdots+B^N\big): \eta \ot \eta \Big]\Big\}\\
 &=\, \la_ 1  \min_{|\eta|=1} \Big\{ \big(B^1+\cdots+B^N\big): \eta \ot \eta  \Big\}\\
&=\, \left(\min_{\ga=1,...,N} \min_{|a|=1} \big\{ A^\ga : a \ot a \big\}\right) \min_{|\eta|=1} \Big\{ \big(B^1+\cdots+B^N\big): \eta \ot \eta \Big\}.
\end{split}
\]
Hence, \eqref{3.7} has been established.

\ms

\noi \textbf{Step 3.} Now we complete the proof of \eqref{3.1} in the general case, by using Steps 1 and 2. We begin by observing the identity
\beq \label{3.8}
\min_{|\eta|=1} \Big\{ \big(B^1+\cdots+B^N\big): \eta \ot \eta \Big\} \, = \, \min_{\ga=1,...,N}\left(\min_{|\eta|=1, \, \eta \in R(B^\ga)} \big\{ B^\ga : \eta \ot \eta \big\}\right).
\eeq
Indeed, \eqref{3.8} follows by applying the Spectral theorem to $B^1,...,B^N$ and by using that by our hypothesis (SH) we have
\[
R(B^1) \oplus \, \cdots \, \oplus R(B^N)\, =\, \R^N, \ \ \ R(B^\ga)\, \bot\, R(B^\de) \ \ \text{ for }\ga \neq \de.
\]
Next, we consider the orthogonal projections on the ranges of $B^\ga$:
\beq \label{3.9}
P^\ga\, :=\, \text{Proj}_{R(B^\ga)},\ \ \ \ga=1,...,N,
\eeq
and we recall that \eqref{3.1} has been established when $N=1$: that is, for any $v \in C^2(\overline{\Om})\cap C^1_0(\Om)$, we have
\beq \label{3.10}
\la_1(A)^2 \int_\Om |D^2v|^2 \, \leq \, \int_\Om |A:D^2v|^2,
\eeq
when $A>0$. Fix now a map $u \in C^2(\overline{\Om})^N\cap C^1_0(\Om)^N $ and apply \eqref{3.10} to $v:=(P^\ga u)_\al$ and $A:=A^\ga$ for indices $\al,\ga$ fixed and then sum with respect to $\al$ and $\ga$, by using that by Step 2 all the $A^\ga$'s have the same first eigenvalue $\la_1$:
\beq \label{3.11}
(\la_1)^2 \int_\Om \sum_\ga \big|D^2 (P^\ga u)\big|^2 \, \leq \, \int_\Om \sum_\ga \big|A^\ga :D^2(P^\ga u)\big|^2.
\eeq
Note now that by perpendicularity we have $P^1 +\cdots +P^N=I$ and hence
\beq \label{3.12}
\big| D^2u\big|^2\, =\, \sum_\ga \big|D^2 (P^\ga u)\big|^2.
\eeq
On the other hand, again by perpendicularity we have $B^\ga B^\de =0$ for $\ga \neq \de$, which implies (we again disengage the summation convention in the equalities right below to avoid confusion)
\[
\begin{split}
\big| \A : D^2u\big|^2\, &=\, \sum_{\al,\ka,\la,i,j,p,q,\ga, \de} \Big( B^\ga_{\al \ka}\, D^2_{ij}u_{\ka} \, A^\ga_{ij} \Big)\, \Big( B^\de_{\al \la}\, D^2_{pq}u_{\la} \, A^\de_{pq}\Big)\\
&=\, \sum_{\al,\ka,\la,i,j,p,q,\ga} \Big( B^\ga_{\al \ka}\, D^2_{ij}u_{\ka} \, A^\ga_{ij}\Big) \, \Big( B^\ga_{\al \la}\, D^2_{pq}u_{\la} \, A^\ga_{pq} \Big) \\
&=\, \sum_\ga \big|B^\ga  D^2 u : A^\ga \big|^2\\
&=\, \sum_\ga \Big|B^\ga P^\ga \big(D^2 u : A^\ga\big)\Big|^2.
\end{split}
\]
Next, for brevity we set
\beq \label{3.13}
\xi^\ga :=\, P^\ga \big(D^2 u : A^\ga\big)\ : \ \ \ \Om\sub \R^n \larrow R(B^\ga), 
\eeq
for $\ga=1,...,N$. Then, by \eqref{3.13}, \eqref{3.12} and Step 2, we may rewrite \eqref{3.11} as
\beq \label{3.11a}
\left(\min_{\ga=1,...,N} \min_{|a|=1} \big\{ A^\ga : a \ot a \big\}\right)^2 \int_\Om \big| D^2u\big|^2\, \leq \, \int_\Om \sum_\ga \big|\xi^\ga|^2.
\eeq
By denoting by ``$\sgn$" the sign function, the above calculation gives (since $\xi^\ga(x) \in R(B^\ga)$ for all $x\in \Om$)
\[
\begin{split}
\big| \A : D^2u\big|^2\, &=\, \sum_\ga \big|B^\ga \xi^\ga \big|^2\\
&=\, \sum_\ga \sup_{\eta \neq 0}\Big( B^\ga : \xi^\ga \ot \frac{\eta}{|\eta|} \Big)^2\\
&\geq \, \sum_\ga  \Big(B^\ga : \xi^\ga \ot \sgn(\xi^\ga) \Big)^2\\
&= \, \sum_\ga  \Big(B^\ga : \sgn(\xi^\ga)\ot \sgn(\xi^\ga) \Big)^2\, |\xi^\ga|^2\\
&\geq \, \sum_\ga  \left[ \min_{\de=1,...,N} \left(\min_{|\eta|=1,\, \eta \in R(B^\de)} \big\{ B^\de : \eta \ot \eta \big\} \right)\right]^2\, |\xi^\ga|^2\\
&= \, \left[\min_{\de=1,...,N} \Big(\min_{|\eta|=1,\, \eta \in R(B^\de)} \big\{ B^\de : \eta \ot \eta \big\} \Big)\right]^2\, \sum_\ga   |\xi^\ga|^2.
\end{split}
\]
By employing the identity \eqref{3.8}, the above estimate gives
\beq \label{3.14}
 \frac{ \big| \A : D^2u\big|^2  }{\, \left[ \min_{|\eta|=1} \Big\{ (B^1+\cdots+B^N): \eta \ot \eta \Big\}\right]^2 \, }\, \geq  \, \sum_\ga   |\xi^\ga|^2.
\eeq
Finally, by \eqref{3.14},  \eqref{3.11a}  and  \eqref{3.7}, the desired estimate \eqref{3.1} follows and Lemma \ref{pr2} ensues.      \qed

\section{Well posedness of the fully nonlinear problem} \label{section4}

We now come to the general fully nonlinear case of the Dirichlet problem \eqref{1.1}. We will utilise the results of Sections \ref{section2} and \ref{section3} plus a result of Campanato on near operators, which is recalled later. Our ellipticity condition of Definition \ref{def2} will work as a ``perturbation device", allowing to establish existence for the nonlinear problem by showing it is ``near" a linear well-posed problem. In view of the well-known problems to pass to limits with weak convergence in nonlinear equations, Campanato's idea furnishes an alternative to the stability problem for nonlinear equations, by avoiding this insuperable difficulty. 

The main result of this paper and this section is the next theorem:

\bt[Existence-Uniqueness for the fully nonlinear problem] \label{th2} Let $\Om\sub \R^n$ be an open, convex, $C^2$ bounded set. Let also $F : \Om \by \R^{Nn^2}_s \larrow \R^N$ a Carath\'eodory map, satisfying Definition \ref{def2} and $F(\cdot,\textbf{0}) \in L^2(\Om)^N$, $n,N\geq 2$. Moreover, suppose that the tensor $\A$ of Definition \ref{def2} satisfies the structural hypothesis (SH).

Then, for any $f\in L^2(\Om)^N$, the problem
\[
\left\{
\begin{array}{l}
F(\cdot,D^2u) \,=\,  f, \ \ \text{ in }\Om,\\
\hspace{36pt} u\,=\, 0, \ \   \text{  on }\p \Om,
\end{array}
\right.
\]
has a unique solution $u$ in the space $(H^2\cap H^1_0)(\Om)^N$, which also satisfies the estimate
\beq \label{4.4a}
\|u\|_{H^2(\Om) }\ \leq \ C\|f\|_{L^2{\Om}},
\eeq
for some $C>0$ depending only on $F$ and $\Om$.  Moreover, for any two maps $w,v \in (H^2\cap H^1_0)(\Om)^N$, we have
\beq \label{5.2}
\|w-v\|_{H^2(\Om)}\, \leq\, C \big\|F(\cdot,D^2w)-F(\cdot,D^2v) \big\|_{L^{2}(\Om)},
\eeq
for some $C>0$ depending only on $F$ and $\Om$. 
\et
 
We note that \eqref{5.2} is a strong uniqueness estimate, which is a form of ``comparison principle in integral norms". Moreover, the restriction to homogeneous boundary condition ``$u=0$ on $\p \Om$" does not harm generality, since the Dirichlet problem we solve is equivalent to a Dirichlet problem with non-homogeneous boundary condition by redefining the nonlinearity $F$: the problem
\[
\left\{
\begin{array}{l}
G(\cdot,D^2u) \,=\,  f, \ \ \text{ in }\Om,\ \ \ \ f\in L^2(\Om)^N\\
\hspace{36pt} u\,=\, g, \ \   \text{  on }\p \Om, \ \ g \in H^2(\Om)^N,
\end{array}
\right.
\]
is equivalent to \eqref{1.1}, by taking $F(x,\X):=G\big(x,\X+D^2g(x)\big)$.

The proof of Theorem \ref{th2} utilises the following result of Campanato taken from \cite{C5}, whose short proof is given for the sake of completeness at the end of the paper in the Appendix:

\bt[Campanato's near operators] \label{th3}  Let $F,A : \mathfrak{X} \larrow X$ be two maps from the set $\mathfrak{X} \neq \emptyset$ to the Banach space $(X,\|\cdot\|)$. Suppose there exists $0<K<1$ such that
\beq \label{5.3}
\Big\| F[u]-F[v]-\big(A[u]-A[v]\big) \Big\| \, \leq\, K \big\| A[u]-A[v] \big\|,
\eeq
for all $u,v \in \mathfrak{X}$. Then, if $A$ is a bijection, $F$ is a bijection as well.
\et

\BPT \ref{th2}. Let $\al$ be the $L^\infty$ function of Definition \ref{def2}. By our assumptions, there exists $C,M>0$ depending only on $F$, such that for any $u\in (H^2\cap H^1_0)(\Om)^N$, we have
\begin{align} \label{5.5}
\big\|\al(\cdot)F(\cdot,D^2u) \big\|_{L^2(\Om)} \, &\leq\, \big\|\al(\cdot) F(\cdot,\textbf{0})\big\|_{L^2(\Om)}\, +\, M\|\al\|_{L^\infty(\Om)} \| D^2u\|_{L^2(\Om)} \\
&=\, \|\al\|_{L^\infty(\Om)} \Big(C\ +\ M \| D^2u\|_{L^2(\Om)} \Big)\nonumber\\
&\leq\, N \Big(1 \, +\, \| u\|_{H^2(\Om)}\Big),\nonumber
\end{align}
for some $N>0$. The last inequality is a consequence of Remark \ref{rem}. Let also $\A \in \R^{Nn \by Nn}_s$ be the tensor given by Definition \ref{def2} corresponding to $F$. Then, we have
\beq \label{5.6}
\|\A:D^2u\|_{L^2(\Om)}\ \leq\ |\A|\, \| D^2u\|_{L^2(\Om)} \ \leq\,  |\A| \| u\|_{H^2(\Om)} .
\eeq
By \eqref{5.5} and \eqref{5.6} we obtain that the operators 
\beq
\left\{
\begin{array}{l}
A[u]\ :=\ \A :D^2u, \ms\\
F[u]\ :=\ \al(\cdot)F(\cdot, D^2u),
\end{array}
\right.
\eeq
map $(H^2\cap H^1_0)(\Om)^N$ into  $L^2(\Om)^N$. Let $u,v \in (H^2\cap H^1_0)(\Om)^N$. By Definition \ref{def2}, we have
\begin{align} 
\Big\|\al(\cdot)\Big(F(\cdot, & D^2u) -  F(\cdot,D^2v)\Big)  -\A: \big(D^2u-D^2v \big)\Big\|_{L^2(\Om)} \nonumber\\
& \leq\ \be \nu(\A) \big\|D^2u-D^2v\big\|_{L^2(\Om)}\ + \ \ga \big\|\A:(D^2u-D^2v)\big\|_{L^2(\Om)}. \nonumber
\end{align}
Since $\A$ satisfies the structural assumption (SH), by the generalised Miranda-Talenti hessian estimate of Lemma \ref{pr2}, we obtain
\begin{align} \label{5.7}
\Big\|\al(\cdot)\Big(F(\cdot, & D^2u) -  F(\cdot,D^2v)\Big)  -\A: \big(D^2u-D^2v \big)\Big\|_{L^2(\Om)} \\
&\leq \ (\be+\ga)\big\|\A:(D^2u -D^2v)\big\|_{L^2(\Om)}. \nonumber
\end{align}
Lemma \ref{le3} implies that the linear operator
\[
A\ :\  (H^2\cap H^1_0)(\Om)^N \larrow  L^2(\Om)^N 
\]
is a bijection. Hence, in view of the inequality \eqref{5.7} and the fact that $\be+\ga<1$, Campanato's Theorem \ref{th3} implies that $F:  (H^2\cap H^1_0)(\Om)^N \larrow  L^2(\Om)^N $ is a bijection as well. As a result, for any $g\in L^2(\R^n)^N$, the PDE system $\al(\cdot)F(\cdot,D^2u)=g$ has a unique solution in $(H^2\cap H^1_0)(\Om)^N$. Since $\al,1/\al \in L^\infty(\Om)$, by selecting $g=\al f$, we conclude that the problem \eqref{1.1} has a unique solution  in $(H^2\cap H^1_0)(\Om)^N$. Finally, by \eqref{5.7} we have
\[
\Big\|F(\cdot,D^2u) -  F(\cdot,D^2v)\Big\|_{L^2(\Om)} \,  \geq \ \frac{1-(\be+\ga)}{\|\al\|_{L^\infty(\Om)}} \, \big\| \A:(D^2u -D^2v) \big\|_{L^2(\Om)}
\]
and by Lemma \ref{pr2} and Remark \ref{rem}, we deduce the estimate
\begin{align} 
\Big\|F(\cdot,D^2u) -  F(\cdot,D^2v)\Big\|_{L^2(\Om)} \, &\geq \ \left(\nu(\A)\frac{1-(\be+\ga)}{\|\al\|_{L^\infty(\Om)}}\right) \big\|D^2u -D^2v\big\|_{L^2(\Om)} \nonumber\\
& \geq\ C\, \| u -v \|_{H^2(\Om)},\nonumber
\end{align}
for some $C>0$. The theorem ensues.      \qed

\section{Stability of the Dirichlet problem} \label{section5}

In this section we discuss an extension of Theorem \ref{th2} in the form of ``stability theorem" for the Dirichlet problem.

\bt[Stability of strong solutions to the Dirchlet problem, cf.\ \cite{K2}] \label{th3} Let $n,N\geq 2$, $\Om\sub \R^n$ a bounded open set and $F,G : \Om \by \R^{Nn^2}_s \larrow \R^N$ Carath\'eodory maps. We suppose that 
\[
F\ :\  (H^2\cap H^1_0)(\Om)^N \larrow  L^2(\Om)^N 
\]
is a bijection. If $G(\cdot,\textbf{0})\in L^2(\Om)^N$ and
\beq \label{6.1}
\underset{x\in \Om}{\ess\,\sup}\sup_{\X\neq\Y}\left| \frac{\big(F(x,\Y)-F(x,\X)\big)-\big(G(x,\Y)-G(x,\X) \big)}{|\Y \, -\, \X|}\right| \ < \ \nu(F)
\eeq
where
\beq \label{6.2}
\nu(F)\ :=\ \inf_{v\neq w}\frac{\big\|F(\cdot,D^2w)-F(\cdot,D^2v) \big\|_{L^{2}(\Om)}}{\|D^2w-D^2v\|_{L^2(\Om)}}\ >\ 0,
\eeq
then, the problem
\[
\left\{
\begin{array}{l}
G(\cdot,D^2u) \,=\,  g, \ \ \text{ in }\Om,\ms\\
\hspace{36pt} u\,=\, 0, \ \   \text{  on }\p \Om,
\end{array}
\right.
\]
has a unique solution $u$ in the space $ (H^2\cap H^1_0)(\Om)^N$, for any given $g\in L^2(\Om)^N$.
\et
 
Theorem \ref{th2} provides sufficient conditions on $F$ and $\Om$ is order to obtain solvability. The theorem says that every $G$ which is ``close to $F$" in the sense of \eqref{6.1}, gives rise to a nonlinear coefficient such that the respective Dirichlet problem is well posed.

\BPT \ref{th3}. 
We denote the right hand side of \eqref{6.1} by $\nu(F,G)$ and we may rewrite  \eqref{6.1}  as
\beq \label{6.3}
0\,<\,  \nu(F,G)\, <\, \nu(F).
\eeq
For any $u,v\in (H^2\cap H^1_0)(\Om)^N$, we have
\begin{align} 
\Big\| F& (\cdot,D^2u) -  F(\cdot,D^2v)  -\big( G(\cdot,D^2u) -  G(\cdot,D^2v)\big) \Big\|_{L^2(\Om)} \nonumber\\
& \leq \left( \underset{\Om}{\ess\,\sup} 
\sup_{\X\neq \Y}\,
\left| \frac{ 
F(\cdot,\Y) -F(\cdot,\X) -
\big( G(\cdot,\Y) -G(\cdot,\X) \big) }{|\Y \, -\, \X|}
\right| \right)
\big\| D^2u-D^2v\big\|_{L^2(\Om)}  \nonumber\\
& = \, \nu(F,G) \big\| D^2u-D^2v\big\|_{L^2(\Om)} \nonumber\\
&\leq\, \frac{\nu(F,G)}{\nu(F)} \big\|F(\cdot,D^2u) -F(\cdot,D^2v)\big\|_{L^2(\Om)}. \nonumber
\end{align}
Hence, we obtain the inequality
\begin{align}  \label{6.4}
\Big\| F (\cdot,D^2u) -  F(\cdot,D^2v)  -&\big( G(\cdot,D^2u) -  G(\cdot,D^2v)\big) \Big\|_{L^2(\Om)} \nonumber\\
&\leq\, \frac{\nu(F,G)}{\nu(F)} \big\|F(\cdot,D^2u) -F(\cdot,D^2v)\big\|_{L^2(\Om)},
\end{align}
which is valid for any $u,v\in (H^2\cap H^1_0)(\Om)^N$. By \eqref{6.1}, Remark \ref{rem} and the inequality above for $v\equiv 0$, we have that $F,G$ map $(H^2\cap H^1_0)(\Om)^N$ into  $L^2(\Om)^N$. By assumption, $F : (H^2\cap H^1_0)(\Om)^N \larrow  L^2(\Om)^N$ is a bijection. Hence, in view of Campanato's Theorem \ref{th3}, inequalities \eqref{6.3} and \eqref{6.4} imply that $G : (H^2\cap H^1_0)(\Om)^N \larrow  L^2(\Om)^N$ is a bijection as well.  The theorem  follows.     \qed

\section*{Appendix: the classical Miranda-Talenti identity and Campanato's nearness theorem}

For completeness and for the convenience of the reader, we include here a variant of the original proof of the Miranda-Talenti inequality (see \cite{M,T}) and the proof of Campanato's theorem on near operators taken from \cite{C5}. 

\bl[Miranda-Talenti identity] \label{MT} Let $\Om \sub \R^n$ be a bounded open set with $C^2$ boundary $\p \Om$. Then, for any $V \in C^1(\overline{\Om})^n$, we have the identity
\[ \tag{A1} \label{A1}
\left(  \text{div} (V) V \, -\, \frac{1}{2}D(|V|^2) \, +\, (n-1)|V|^2 \H\right) \cdot N\, =\, 0,\ \ \text{ on }\p\Om,
\]
where $N$ is the unit outward pointing normal vector field of $\p\Om$, $\H$ is the mean curvature vector of $\p\Om$ and $V$ satisfies
\[
N^\bot V\,=\, 0
\]
where $N^\bot=I-N\ot N$. Moreover, for any $v \in (H^2\cap H^1_0)(\Om)$, the identity \eqref{3.2} holds true.
\el

\BPL \ref{MT}. The essential part of the result is identity \eqref{A1}. Assuming it, let us show how we can obtain \eqref{3.2}: for any $u\in C^\infty(\overline{\Om})$ with $u=0$ on $\p\Om$, we have the elementary identity
\[
 | \De u |^2\, -\,  |D^2u  |^2\,=\, \text{div}(\De u Du) \, -\, \frac{1}{2}\De (|Du|^2)
\]
which by integration and Gauss' theorem  gives
\[
\int_\Om\Big(  | \De u |^2\, -\,  |D^2u  |^2 \Big)d\mL^n\, =\, \int_{\p\Om}  \left( \De u Du  \, -\, \frac{1}{2}D (|Du|^2) \right) \cdot N \, d\mH^{n-1}.
\]
Then, by applying \eqref{A1} to $V=Du$ and by using that $N^\bot Du=0$ which is a consequence of the vanishing on $u$ on $\p\Om$, we obtain \eqref{3.2}.

Let us now demonstrate \eqref{A1}. For the rest of the paper, the indices $i,j,k,...$ run in $\{1,...,n-1\}$ and the $\al,\be,\ga$ in $\{1,...,n\}$. Let 
\[
X\  :\  T\sub \R^{n-1}\larrow X(T) \sub \p\Om 
\]
be a local parametrisation of the boundary. Then, for any $x=X(t) \in X(T)$, $DX(t)=D_iX_\al e^\al \ot e^i$ spans the tangent space $T_x\p\Om$. For any fixed $f\in C^1(\overline{\Om})$, we now derive an identity which relates the Euclidean gradient $Df$ of (an extension on $\R^n$ of) $f$  to the geometry of the domain. Let $g$ be the induced Riemann metric on $\p\Om$ and $\nabla f$ the Riemannian gradient on $\p\Om$. Then, $g_{ij}(t)=D_i X_\al(t) D_j X_\al(t)$  and if $g^{ij}$ denotes the components of the inverse $g^{-1}$, we have
\[
\nabla f(x)\, =\, D_\al f(x)\, D_i X_\al(t)\, g^{ij}(t)\, D_jX(t),
\]
at $x=X(t) \in X(T) \sub \p\Om$. Moreover, since $Df$ and $\nabla f$ are also satisfy
\[
\nabla f\,=\, N^\bot Df\, =\, Df\, -\, N\ot N Df,
\]
we obtain 
\[ \tag{A2} \label{A2}
D_\al f(x)\, =\, N_\al(x) \, N_\be (x)\,D_\be f(x)\ +\  D_\al f(x)\, D_i X_\al(t)\, g^{ij}(t)\, D_jX(t),
\]
for any $x=X(t) \in X(T) \sub \p\Om$. We recall now that the 2nd fundamental form on $\p\Om$ is given in local coordinates by
\[
b_{ij}\circ X \, =\, -D_i (N_\al \circ X) \, D_jX_\al  
\]
and the mean curvature (with respect to the orientation defined by the outward pointing nornal $N$) by
\[
H\circ X\, =\, \frac{1}{n-1}g^{ij}b_{ij}\, =\, -\frac{1}{n-1}g^{ij}\, D_i (N_\al \circ X) \, D_jX_\al .
\]
In order to prove \eqref{A1}, for $V$ as in the statement of the lemma and in view of \eqref{A2}, on $X(T) \sub \p\Om$ we have
\[
\begin{split}
\text{div}(V) V\cdot N \, -\,   \frac{1}{2} & D(|V|^2) \cdot N\  = \, -\, N_\al V_\be\, (D_\al V_\be)  \  + \  V_\be \,N_\be\, N_\al \, N_\ga \, (D_\ga V_\al) \ \\
&\ +\, V_\be\, N_\be\, (D_\ga V_\al)\, \big( D_jX_\ga \circ X^{-1}\big)\big( g^{ij}\circ X^{-1}\big) \, \big( D_iX_\al \circ X^{-1} \big) .
\end{split}
\]
By using that $V_\al =N_\al N_\be V_\be$, the first two terms cancel and by further using that $V_\al =|V|N_\al $, we get
\[
\begin{split}
\text{div} (V)  V & \cdot  N \, -\,   \frac{1}{2}  D(|V|^2) \cdot N \\
& =\, V_\be\, N_\be\, (D_\ga V_\al)\, \big( ( D_jX_\ga) \circ X^{-1}\big)\big( g^{ij}\circ X^{-1}\big) \, \big( (D_iX_\al ) \circ X^{-1} \big) \\
&=\, |V|\, \big( g^{ij}\circ X^{-1}\big) \, \Big( D_j (V_\al \circ X) \circ X^{-1} \Big)  \, \big( D_iX_\al \circ X^{-1} \big)\\
&=\, |V|\, \big( g^{ij}\circ X^{-1}\big) \, \Big( D_j \Big[ | V \circ X|(N_\al \circ X) \Big] \circ X^{-1} \Big)  \, \big( D_iX_\al \circ X^{-1} \big),
\end{split}
\]
on $X(T) \sub \p\Om$. By expanding the derivative and using that $N \cdot (D_k X)\circ X^{-1}=0$, we have
\[
\begin{split}
\text{div}  (V)  V\cdot  N \, -\,   \frac{1}{2}  &D(|V|^2) \cdot N \\
&=\, |V|\, \big( g^{ij}\circ X^{-1}\big) \, \Big( 
 | V| \big(D_j(N_\al \circ X)\big)\circ X^{-1}  \\
 &\ \ \ \ \ +\  N_\al\, D_j\big(| V\circ X|\big)\circ X^{-1} 
\Big)  
\, \big( (D_iX_\al ) \circ X^{-1} \big)\\
&=\, |V|^2\, \big( g^{ij}\circ X^{-1}\big)  \, \big( (D_iX_\al) \circ X^{-1} \big) \, D_j(N_\al \circ X) \circ X^{-1} \\
&=\, -\,|V|^2\, \big(g^{ij}\circ X^{-1}\big)\,   \big( b_{ij}\circ X^{-1}\big)\\
&=\, -\,(n-1)|V|^2H\\
&=\, -\,(n-1)|V|^2\H \cdot N,
\end{split}
\]
on $X(T) \sub \p\Om$, since $H=\H \cdot N$. The conclusion follows by covering $\p\Om$ by coordinate neighbourhoods.             \qed

\ms

We now give the proof of Campanato's theorem.

\BPT \ref{th3}. It suffices to show that for any $f\in X$, there is a unique $u\in \mathfrak{X}$ such that
\[
F[u]\, =\, f.
\]
In order to prove that, we first turn $\mathfrak{X}$ into a complete metric space, by pulling back the structure from $X$ via $A$: for, we define the distance
\[
d(u,v)\, :=\, \big\| A[u]-A[v]\big\|.
\]
Next, we fix an $f\in X$ and define the map
\[
T\ : \ \mathfrak{X} \larrow \mathfrak{X}\ , \ \ \ T[u]\, :=\, A^{-1}\Big(A[u]-\big(F[u]-f \big) \Big).
\]
We conclude by showing that $T$ is a contraction on $(\mathfrak{X},d)$, and hence has a unique $u\in \mathfrak{X}$ such that $T[u]=u$. The latter equality is equivalent to $F[u]=f$, and then we will be done.  Indeed, we have that
\begin{align}
d\Big( T[u],T[v] \Big)  \, &=\ \Big\|  \left(A[u]-\big(F[u]-f \big) \right) \, -\,  \left(A[v]-\big(F[v]-f \big) \right) \Big\| \nonumber \\
 &=\ \Big\| A[u] -A[v] -\big(F[u]  - F[v]\big) \Big\|, \nonumber 
\end{align}
and hence
\begin{align}
d\Big( T[u],T[v] \Big)  \ \, 
&\!\!\!\!\overset{\eqref{5.3}}{\leq}   K \big\| A[u]-A[u] \big\| \nonumber\\
  &= \  K \, d(u,v). \nonumber
\end{align}
Since $K<1$, the conclusion follows and the theorem ensues.              \qed

\ms

\ms

\bibliographystyle{amsplain}

\end{document}